# AZB Rectangle Shrinkage Method and Heterogeneous Computing Accelerated Full Image Theory Method Ray Tracing Enabling Complex and Massive Outdoor 6G Propagation Modeling

Yongwan Kim, *Member, IEEE*, Hyunjun Yang, *Member, IEEE,* and Jungsuek Oh, *Senior Member, IEEE*

*Abstract*—Until now, despite their high accuracy, the utilization of the conventional image theory method ray tracers was limited to simple simulation environments with small number of field observation points and low maximum ray bouncing order due to their poor computational efficiency. This study presents a novel full-3D AZB rectangle shrinkage method and heterogeneous computing accelerated image theory method ray tracing framework for complex and massive outdoor propagation modeling. The proposed framework is divided into three parts: 1. Visibility preprocessing part. 2. Visibility tree generation part: in this part, a novel AZB rectangle shrinkage method that accelerates and reduces generation speed and size of visibility tree is proposed. 3. Shadow testing and field calculation part: in this part, a heterogeneous computing algorithm that can make possible to handle a large amount of field observation points is proposed. It is demonstrated that the proposed framework is faster more than 651 times than the image theory method solver of WinProp. Also, it is confirmed that the proposed ray tracing framework can handle 1km x 1km wide and dense urban outdoor simulation with up to the maximum ray bouncing order of 6 and thousands of field observation points. The proposed ray tracing framework would be a cornerstone of future image theory method ray tracing techniques for complex and massive scenarios that was exclusive to the shooting and bouncing rays method ray tracers.

## I. INTRODUCTION

MILIMETER-WAVE (MMWAVE) COMMUNICATION SYSTEMS requiring the sophisticated manipulation of LOS and NLOS propagation paths channel essentially demand for the characterization of deterministic electromagnetic channel is based on channel data obtained from not only measurement campaigns but also ray tracing (RT) electromagnetic (EM) analysis [1]-[5]. Especially, the RT is extensively utilized over the mmWave frequencies and would be a preferable tool for following 6G channel modeling and cell coverage analysis due to its convenience, while actual channel measurement on 6G frequency band is very complex and time-consuming [6]. For outdoor channel modeling and cell coverage analysis using RT, EM analysis for a large amount of field observation point (FOP) on complex and massive simulation environment is performed [1]-[5].

There are two representative methods of RT. First one is shooting and bouncing rays (SBR) method, the most widely used ray tracing method that generally provides high computational efficiency; however, the accuracy is limited owing to phase error and cutoff of the ray tube [7]. The second one is image theory (IT) method which presents high accuracy but a relatively low computational efficiency, compared with the SBR method owing to its shadow testing (ST) and visibility tree generation algorithm [7].

Recently, frequency range we use for communication systems is getting higher and higher to achieve broad bandwidth, higher data rate and low latency, etc. If an analysis frequency increases, the SBR method generally produces higher phase error which is a by-product of the ray tube concept [8]-[10]. To maintain accuracy of SBR method when analysis frequency increases, it should launch more rays or apply additional optimization procedure such as stationary point search [8], [9], resulting in larger computation task and complication of algorithm. Also, the SBR method has ray tube cutoff error, which occurs when the size of wavefront of the ray tube becomes larger than resolution of simulation environment [11]. As operating frequency increases, antenna main beam width narrows to overcome large path loss [6], so that the number of important strong rays decreases a lot relative to weak rays. If the ray tube cutoff error occurs for these few important rays, the overall analysis results would be invalid. Therefore, the total number of source rays should increase to achieve acceptable accuracy, leading to larger computational load [11]. Furthermore, at this point, there is an ambiguity for determining how many rays should be launched, to escape ray tube cutoff

Manuscript received August xx, 2023; This work was supported by the Institute of Information & Communications Technology Planning & Evaluation (IITP) grant funded by the South Korea government (MIST) (No.2019-0-00098, Advanced and Integrated Software Development for Electromagnetic Analysis) *(Corresponding author: Jungsuek Oh.)*
The authors are with the Institute of New Media and Communications (INMC) and the Department of Electrical and Computer Engineering, Seoul National University, Seoul 08826, South Korea (e-mail: yongwankim@snu.ac.kr, jungsuek@snu.ac.kr)

 

and achieve sufficient accuracy for various analysis frequencies and simulation environments.

On the other hand, for the IT method, the accuracy is affected by only geometric mesh approximations and floating-point number precision of computers [11]. Therefore, it presents zero systemic error without increase of computational load and complication of algorithm even if the analysis frequency increases. Also, its accuracy is independent with the resolution of simulation environment and antenna main beam width, so that explicit simulation is possible without any ambiguity, such as ray tube cutoff and how many rays should be launched to achieve sufficient accuracy [11]. Considering facts above, it can be seen that the IT method is a promising RT method for the analysis of the 6G communication systems. However, the conventional IT method suffers from prohibitively long computation time when the propagation scenario treats high order of maximum ray bounce, numerous FOPs and facets modeling simulation environment, so that its utilization is limited to simple simulation environments so far [10], [12].

To overcome this long computational time, Catedra et al. [13] utilized the angular Z-buffer (AZB) to the IT method and achieved computation time reduction of 90% [10]. The AZB combined IT method divides the space into multiple angular space around a source, i.e. Tx or image Tx. After that, the source radiates a beam to the AZB rectangle enclosing the reflection facet. And the same procedure is proceeded for facets inside the beam up to the maximum ray bouncing order. At this moment, the conventional AZB method creates AZB rectangle enclosing the entire reflecting facet including the portion that is not illuminated by the beam of former bouncing order, i.e. portion of the facet that does not participate the multiple reflection. As a result, a large amount of unnecessary computational task and memory usage, which increases exponentially as the ray bouncing order, are produced. If a simulation scenario includes diffraction, this inefficiency much more worsen because diffraction beam radiated to the AZB rectangle generally illuminates much more facets than reflection beam, resulting in infeasibility of massive analysis. Therefore, to avoid this inefficiency, the beam should be radiated only to the portion that is illuminated by the beam of former bouncing order. To make it possible, we propose an AZB rectangle shrinkage (ARS) method as the first contribution of this paper.

The conventional AZB method creates DAZB matrices storing information about which facets are contained in which angular region. Using DAZB matrices, the ST can be accelerated by considering only facets inside the angular region where ray belongs. However, because the DAZB matrices should be created for every Tx and image Tx, when a high order of ray bouncing order is considered, it can occupy prohibitively large amount of memory, precluding an acceleration through parallel computing. In this context, we propose a novel heterogeneous computing, i.e., CPU and GPU parallel computing, based ARS-accelerated IT (HAIT) method RT framework which uses the AZB rectangle concept but does not use the DAZB matrix, as the second contribution of this paper. The HAIT has an additional advantage that it can handle a large amount of FOPs using GPU parallel computing. Although parallel computing approach has been widely used to resolve the limited computation speed of RT recently, in most cases, it has been applied to the SBR method [14]-[24]. Conversely, for the IT method, only few studies on parallel computing acceleration exist and these have obvious limitations. For example, in [10], the GPU-based kD-tree-accelerated beam tracing (GKBT) method is presented. Although the GKBT method was approximately ten times faster than ray tracing techniques before the GKBT was developed, it cannot support the reflection of diffracted fields owing to the use of the ray congruence concept for diffraction. The reflection of diffracted fields is crucial in outdoor wave propagation; hence, the GKBT is unsuitable for analyzing outdoor environments. On the other hand, our technique supports diffraction-reflection analysis with maximally 6 reflection and 1 diffraction, so that it is feasible for outdoor scenario analysis. Moreover, although [25] and [26] present the GPU-accelerated IT method RT technique, these studies mainly focus on visibility preprocessing and the diffuse scattering effect, respectively, and only limited information on the RT acceleration technique is provided. Considering facts above, the second contribution of this paper would be a cornerstone of the future IT method RT technique for massive analysis.

The paper is organized as follows. In Section II, the numerical formulations of ARS method for multiple reflection, reflection-diffraction, and diffraction-reflection propagation sequences are introduced. In Section III, the HAIT method RT framework is elaborated in detail. In Section IV, the accuracy and efficiency of ARS and HAIT method are validated using simple and complex urban outdoor scenario. In Section V, a conclusion is described.

## II. AZB Rectangle Shrinkage (ARS) Method

In this section, we propose the ARS method which reduces the size of the AZB rectangle as ray bouncing order increases, leading to exponential reduction of size and acceleration of generation speed of the visibility tree as maximum bouncing order increases, compared to the conventional AZB IT method. The ARS consists of three major cases: *A. Multiple reflection. B. Reflection-diffraction. C. Diffraction-reflection.*

We note that the general spherical coordinate system in [10] is used for ARS.

### A. Multiple reflection

A simple 2D description of second order multiple reflection of ARS is shown in Fig. 1 where first and second reflecting facets are represented by 1st and 2nd facet, respectively. Firstly, for first order reflection, we generate an image Tx for a 1st facet which is represented as 1st image Tx in Fig. 1. Subsequently, from the image Tx, we see the 1st facet and determine its $\phi - \theta$ AZB rectangle, which is a rectangle determined by four spherical angular margins ($\phi_{min}, \phi_{max}, \theta_{min},$ and $\theta_{max}$) of the reflecting facet (see Fig. 2). This AZB rectangle and image Tx form the 1st beam as shown in Fig. 1. Next, we determine the AZB rectangle of the facet inside the 1st beam, i.e., 2nd facet,

  3

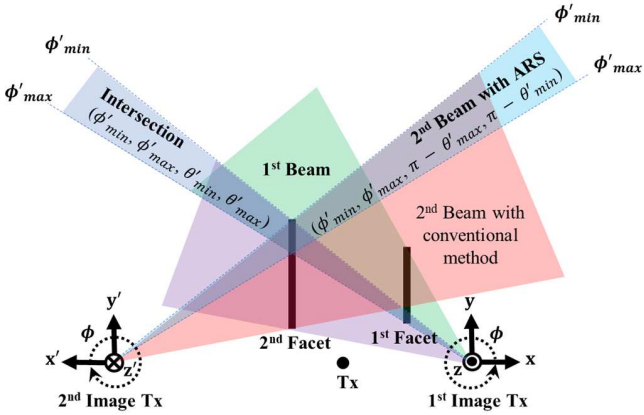

Fig. 1. 2D description of second order multiple reflection of ARS.

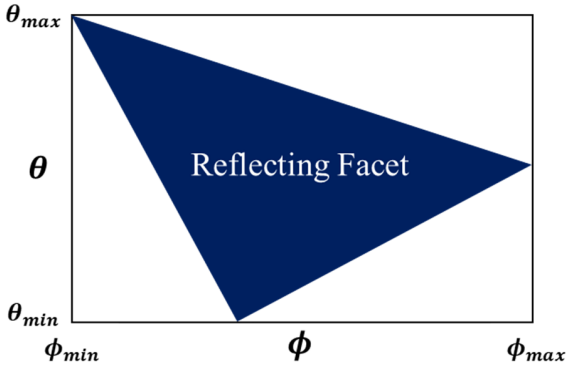

Fig. 2. Example of $\phi - \theta$ AZB rectangle.

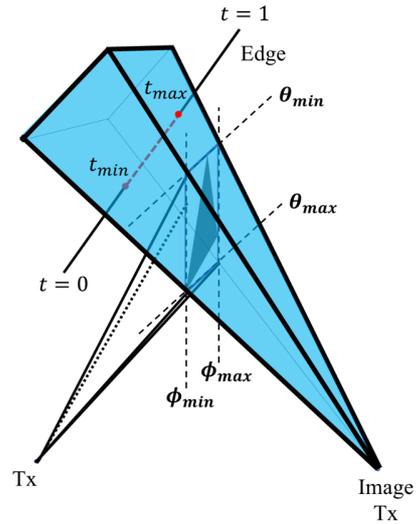

Fig. 3. Problem description of reflection-diffraction case of ARS

and find the intersection between two AZB rectangles, i.e., $\phi'_{min}, \phi'_{max}, \theta'_{min}$, and $\theta'_{max}$ in Fig. 1. Subsequently, to analyze second order reflection, we generate another image Tx about the 2$^{nd}$ facet and transform basis as 2$^{nd}$ image Tx shown in Fig. 1. The first and second new basis, i.e., $\mathbf{x}'$ and $\mathbf{y}'$, are $\mathbf{x}$ and $\mathbf{y}$ mirrored to the 2$^{nd}$ facet, respectively, and the last new basis, i.e., $\mathbf{z}'$, equals $\mathbf{x}' \times \mathbf{y}'$. Then, we can directly use the intersected AZB rectangle, that is, reduced AZB rectangle, for the second bouncing order after changing $\theta'_{min}$ and $\theta'_{max}$ to $\pi - \theta'_{max}$ and $\pi - \theta'_{min}$, respectively, which means that 2$^{nd}$ beam is launched from 2$^{nd}$ image Tx to the AZB rectangle having angular margins of $\phi'_{min}, \phi'_{max}, \pi - \theta'_{max},$ and $\pi - \theta'_{min}$ as shown in Fig. 1. Subsequently, an identical procedure is performed up to the maximum ray bouncing order.

From Fig. 1, it can be seen that size of 2$^{nd}$ beam with ARS is significantly reduced compared to the one with the conventional method and this efficiency improvement would increase exponentially as the bouncing order increases. Also, although only 2D case is illustrated in Fig. 1 for simplicity, this ARS method for multiple reflection can be applied to facets with any 3D locations and orientations.

*B. Reflection-diffraction*

A simple problem description of reflection-diffraction case of ARS is shown in Fig. 3. In this figure, an AZB rectangle is formed for a facet and a reflection beam, i.e., blue beam, is launched from an image Tx to the AZB rectangle. At this time, an edge penetrates the reflection beam and the portion of the edge participating the reflection-diffraction is only the portion inside the reflection beam, which is represented as a red dotted line in Fig. 3. This is the portion between $t_{min}$ and $t_{max}$ where $t$ is the normalized length of the edge, which means the start point and end point of the edge is represented as $t = 0$ and $t = 1$, respectively, as shown in Fig. 3. Although the exact portion of the edge that actually participates the diffraction can be easily found using simple beam tracing algorithm in [10] when a first reflection-diffraction is considered, the beam tracing algorithm requires complex beam splitting procedure when multiple reflection-diffraction is considered, resulting in exponential increase of computational load and algorithm complexity. Instead, we use the AZB rectangle concept to reflection-diffraction and we can directly use the shrunk AZB rectangle through multiple reflections described in Section II-A for multiple reflection-diffraction. The $t_{min}$ and $t_{max}$ can be expressed as follows:

$$t_{min} = \max(t_{\phi min}, t_{\theta min}, 0) \quad (1)$$
$$t_{max} = \min(t_{\phi max}, t_{\theta max}, 1) \quad (2)$$

where $t_\phi$ and $t_\theta$ are $t$ values where the edge intersects $\phi$ and $\theta$-axis sides of the reflection AZB rectangle, respectively.

Firstly, to formulate $t_\phi$, we represent the edge line as $\mathbf{E}$ in (3).

$$\mathbf{E} = \mathbf{P_1} + t(\mathbf{P_2} - \mathbf{P_1}) \quad (3)$$

Then, if we represent the location of the image Tx as $\mathbf{S}$, a vector from the image Tx to a point on the edge line is represented as follows:

$$\mathbf{SE} = \mathbf{A} + t\mathbf{B} \quad (4)$$

where $\mathbf{A}$ and $\mathbf{B}$ are $\mathbf{P_1} - \mathbf{S}$ and $\mathbf{P_2} - \mathbf{P_1}$, respectively. Now, we formulate $\tan \phi$ of a point on the edge line considering the image Tx as the origin as follows:

$$\tan \phi = \frac{A_y + tB_y}{A_x + tB_x} \quad (5)$$

After performing some algebra with (5), we get

$$t(\phi) = \frac{-A_x \tan \phi + A_y}{B_x \tan \phi - B_y} \quad (6)$$

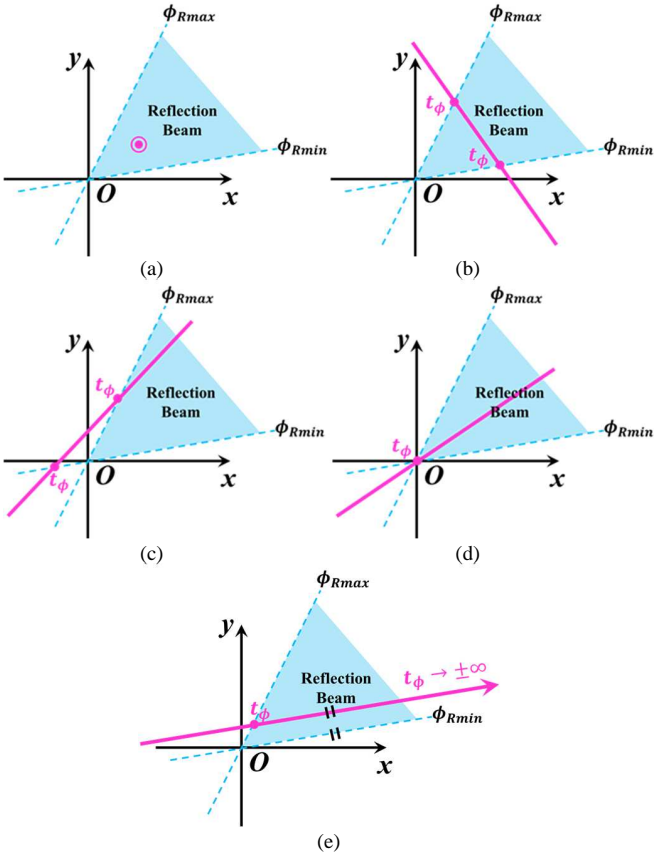

Fig. 4. 5 possible cases where whole or part of edge line can be located inside the $\phi$ domain of the reflection AZB rectangle.

At this moment, we note that $t(\phi)$ equals $t(\phi + \pi)$ because $\tan\phi$ equals $\tan(\phi + \pi)$; therefore, a validation procedure of the $t$ value should be performed by confirming atan2 function of $\mathbf{A} + t(\phi)\mathbf{B}$ equals $\phi$. There are 5 possible cases where whole or a part of edge line can be located inside the $\phi$ domain of the reflection AZB rectangle, and these are shown in Fig. 4. Fig. 4(a) represents a case where the edge line is directed to the z-axis, and whole edge line is inside the $\phi$ domain of the AZB rectangle, resulting in $t_{\phi min} \to -\infty$ and $t_{\phi max} \to \infty$. Fig. 4(b) shows the case where the edge line penetrates the $\phi$ domain of the AZB rectangle, so that two $t_\phi$, i.e., $t(\phi_{Rmin})$ and $t(\phi_{Rmax})$ where $\phi_{Rmin}$ and $\phi_{Rmax}$ are $\phi$-domain angular margins of AZB rectangle, are assigned to the $t_{\phi min}$ and $t_{\phi max}$, respectively, depending on their magnitude. Fig. 4(c) shows the case where one $t_\phi$ occurs at one $\phi$-domain angular margin and the other $t_\phi$ occurs at the other angular margin $+ \pi$, so that the first $t_\phi$ and $\pm\infty$ are allocated to the $t_{\phi min}$ and $t_{\phi max}$, respectively, depending on their magnitude. Fig. 4(d) shows the case where projection of the edge to the $xy$-plane ($z = 0$) passes through the origin, so that $t_\phi = t(\phi_{Rmin}) = t(\phi_{Rmax})$ and this $t_\phi$ and $\pm\infty$ is allocated to $t_{\phi min}$ and $t_{\phi max}$, respectively, depending on their magnitude. Fig. 4(c) shows the case where one $t_\phi$ occurs at one $\phi$-domain angular margin and the other angular margin is parallel with the edge line, so that the $t_\phi$ and $\pm\infty$ is allocated to $t_{\phi min}$ and $t_{\phi max}$, respectively, depending on their magnitude. The pseudocode for the detailed determination procedure of $t_{\phi min}$ and $t_{\phi max}$ is shown in Algorithm 1. We note that if the reflecting facet is located at $\pm z$

---

**Algorithm 1** Determination of $t_{\phi min}$ and $t_{\phi max}$

**Input:**
   $\phi_{Rmin}, \phi_{Rmax}, \mathbf{P_1}, \mathbf{P_2},$ and $\mathbf{S}$
**Output:**
   $t_{\phi min}$ and $t_{\phi max}$

1: **function** tPhiMinMax
2:   **if** $P_{1x} - P_{2x} = P_{1y} - P_{2y} = 0$ **then**
3:     **if** atan2$(P_{1y} - S_y, P_{1x} - S_x)$ is between $\phi_{Rmin}$ and $\phi_{Rmax}$ **then**
4:       $t_{\phi min}, t_{\phi max} \leftarrow -\infty, \infty$     ▷ Fig. 4(a)
5:     **else**
6:       $t_{\phi min}, t_{\phi max} \leftarrow$ NaN, NaN
        ▷ Edge line is located outside the reflection beam
7:     **end if**
8:     **return** $t_{\phi min}, t_{\phi max}$
9:   **end if**
10:  $\mathbf{A}, \mathbf{B} \leftarrow \mathbf{P_1} - \mathbf{S}, \mathbf{P_2} - \mathbf{P_1}$
11:  $t_{\phi 1}, t_{\phi 2} \leftarrow t(\phi_{Rmin}), t(\phi_{Rmax})$
12:  $\mathbf{d_1}, \mathbf{d_2} \leftarrow \mathbf{A} + t_{\phi 1}\mathbf{B}, \mathbf{A} + t_{\phi 2}\mathbf{B}$
13:  $\phi_1, \phi_2 \leftarrow$ atan2$(d_{1y}, d_{1x})$, atan2$(d_{2y}, d_{2x})$
14:  **if** $\phi_1 = \phi_{Rmin}$ and $\phi_2 = \phi_{Rmax}$ **then**   ▷ Fig. 4(b)
15:    $t_{\phi min}, t_{\phi max} \leftarrow \min(t_{\phi 1}, t_{\phi 2}), \max(t_{\phi 1}, t_{\phi 2})$
16:  **else if** $(\phi_1 = \phi_{Rmin}$ and $\phi_2 = \phi_{Rmax} + \pi)$ or $(\phi_1 = \phi_{Rmin} + \pi$ and $\phi_2 = \phi_{Rmax})$ **then**   ▷ Fig. 4(c)
17:    $t_{\phi proper} \leftarrow$ The $t$ value that matches $\phi_{Rmin}$ or $\phi_{Rmax}$ between $t_{\phi 1}$ and $t_{\phi 2}$
18:    $t_{\phi improper} \leftarrow$ The $t$ value that matches $\phi_{Rmin} + \pi$ or $\phi_{Rmax} + \pi$ between $t_{\phi 1}$ and $t_{\phi 2}$
19:    **if** $t_{\phi proper} - t_{\phi improper} > 0$ **then**
20:      $t_{\phi min}, t_{\phi max} \leftarrow t_{\phi proper}, \infty$
21:    **else**
22:      $t_{\phi min}, t_{\phi max} \leftarrow -\infty, t_{\phi proper}$
23:    **end if**
24:  **else if** $t_{\phi 1} = t_{\phi 2}$ **then**
25:    $\mathbf{C} \leftarrow \mathbf{A} + (t_{\phi 1} + \alpha)\mathbf{B}$   ▷ $\alpha$ is arbitrary positive constant
26:    $\phi_0 \leftarrow$ atan2$(C_y, C_x)$
27:    **if** $\phi_0$ is between $\phi_{Rmin}$ and $\phi_{Rmax}$ **then**   ▷ Fig. 4(d)
28:      $t_{\phi min}, t_{\phi max} \leftarrow t_{\phi 1}, \infty$
29:    **else if** $\phi_0$ is between $\phi_{Rmin} + \pi$ and $\phi_{Rmax} + \pi$ **then**
      ▷ Fig. 4(d)
30:      $t_{\phi min}, t_{\phi max} \leftarrow -\infty, t_{\phi 1}$
31:    **else**
32:      $t_{\phi min}, t_{\phi max} \leftarrow$ NaN, NaN
33:    **end if**
34:  **else if** $(\phi_1 = \phi_{Rmin}$ and $t_{\phi 2} \to \pm\infty)$ or $(\phi_2 = \phi_{Rmax}$ and $t_{\phi 1} \to \pm\infty)$   ▷ Fig. 4(e)
35:    $t_{\phi finite} \leftarrow$ The finite $t$ value between $t_{\phi 1}$ and $t_{\phi 2}$
36:    $\mathbf{C} \leftarrow \mathbf{A} + (t_{\phi finite} + \alpha)\mathbf{B}$
      ▷ $\alpha$ is arbitrary positive constant
37:    **if** atan2$(C_y, C_x)$ is between $\phi_{Rmin}$ and $\phi_{Rmax}$ **then**
38:      $t_{\phi min}, t_{\phi max} \leftarrow t_{\phi finite}, \infty$
39:    **else**
40:      $t_{\phi min}, t_{\phi max} \leftarrow -\infty, t_{\phi finite}$
41:    **end if**
42:  **else**
43:    $t_{\phi min}, t_{\phi max} \leftarrow$ NaN, NaN
44:  **end if**
45:  **return** $t_{\phi min}, t_{\phi max}$
46: **end function**

direction from the image Tx, i.e., a ray with direction of $\pm z$ launched from image Tx intersects the reflecting facet, we do not conduct Algorithm 1 and just set $t_{\phi min}$ and $t_{\phi max}$ to $-\infty$ and $\infty$, respectively, because $\phi_{Rmin}$ and $\phi_{Rmax}$ are 0 and $2\pi$, respectively, at this case.

Next, to derive $t_\theta$, we formulate $\cot\theta$ of a point on the edge line considering the image Tx as the origin as follows:

$$\cot\theta = \frac{A_z + tB_z}{\sqrt{(A_x + tB_x)^2 + (A_y + tB_y)^2}} \quad (7)$$

After performing some algebra with (7), we get

$$t(\theta) = \frac{-\beta \pm \sqrt{\beta^2 - 4\alpha\gamma}}{2\alpha} \quad (8)$$

where

$$\alpha = (B_x^2 + B_y^2)\cot^2\theta - B_z^2 \quad (9)$$
$$\beta = 2[(A_xB_x + A_yB_y)\cot^2\theta - A_zB_z] \quad (10)$$
$$\gamma = (A_x^2 + A_y^2)\cot^2\theta - A_z^2 \quad (11)$$

It should be noted that $t(\theta)$ equals $t(\pi - \theta)$ because $\cot^2\theta$ equals $\cot^2(\pi - \theta)$; therefore, a selection procedure of correct $t$ value should be performed by substituting the calculated two $t$ to the following equation.

$$\theta(t) = \cot^{-1}\frac{A_z + tB_z}{\sqrt{(A_x + tB_x)^2 + (A_y + tB_y)^2}} \quad (12)$$

There are three mathematical facts: 1. $t(\theta) = t(\pi - \theta)$. 2. $t(\theta)$ has maximally two solutions. 3. The limit of (12) as $t$ approaches to infinity is expressed as follows:

$$\lim_{t\to\pm\infty}\theta(t) = \cot^{-1}\left(\pm B_z/\sqrt{B_x^2 + B_y^2}\right) \quad (13)$$

Based on mathematical facts above, we can derive other two mathematical facts: 1. $\theta(t)$ can have maximally one extreme value. 2. Local minimum and local maximum value of $\theta(t)$ cannot exist at the area where $\theta(t) > \pi/2$ and $\theta(t) < \pi/2$, respectively. Based on these two mathematical facts, we derive total 20 possible cases, which is shown in Fig. 5, where whole or a part of edge line can be located inside the $\theta$ domain of the reflection AZB rectangle.

The pseudocode for the detailed determination procedure of $t_{\theta min}$ and $t_{\theta max}$ is shown in Algorithm 2. We note that there occur beam splitting at some cases, e.g., Fig. 5(h), (i); however, for simplicity, we ignored it and merged the beams in this paper. Applying the beam splitting would improve the performance of algorithm and it would be a promising future work.

*C. Diffraction-reflection*

To conveniently treat diffraction $\phi - t$ AZB rectangle, we transform basis as shown in Fig. 6. For the first diffraction $\phi - t$ AZB rectangle, $\phi_{min} = 0$, $\phi_{max} = n\pi$, and $t_{min} = 0$ and $t_{max} = 1$. Subsequently, we determine the facet inside the $\phi - t$ AZB rectangle. Next, we determine the second $\phi - t$ AZB rectangle for the facet seen from the edge. $\phi_{min}$ and $\phi_{max}$ of the second AZB rectangle can be easily determined by calculating $\phi$ of three vertices of the facet and comparing their magnitude. However, to determine the exact values of $t_{min}$ and $t_{max}$ of the mesh, nonlinear equation should be solved and it is

---

**Algorithm 2** Determination of $t_{\theta min}$ and $t_{\theta max}$

**Input:**
  $\theta_{Rmin}, \theta_{Rmax}$, **A**, and **B**
**Output:**
  $t_{\theta min}$ and $t_{\theta max}$

1: **function** tThetaMinMax
2:   $t_{Rmin+}, t_{Rmin-} \leftarrow t(\theta_{Rmin})$ with $\pm$ sign
3:   $t_{Rmax+}, t_{Rmax-} \leftarrow t(\theta_{Rmax})$ with $\pm$ sign
4:   **if** all $t_{Rmin\pm}$ and $t_{Rmax\pm}$ invalid **then**
5:     **if** $\mathbf{A} \times \mathbf{B} \neq \mathbf{0}$ **then**    ▷ Expanded edge line not intersects origin
6:       **if** $\theta_{Rmin} < \theta(\delta) < \theta_{Rmax}$ **then**
        ▷ Fig. 5(a), $\delta$: arbitrary real value
7:         $t_{\theta min}, t_{\theta max} \leftarrow -\infty, \infty$
8:       **else**
9:         $t_{\theta min}, t_{\theta max} \leftarrow$ NaN, NaN
10:       **end if**
11:     **else** ▷ Expanded edge line intersects origin
12:       $t_0 \leftarrow value\ that\ satisfies\ \mathbf{A} + t_0\mathbf{B} = 0$
13:       **if** $\theta_{Rmin} < [\theta(t_0 + \varepsilon)\ \text{and}\ \theta(t_0 - \varepsilon)] < \theta_{Rmax}$ **then**
        ▷ Fig. 5(b), $\varepsilon$: arbitrary positive value
14:         $t_{\theta min}, t_{\theta max} \leftarrow -\infty, \infty$
15:       **else if** $\theta_{Rmin} < \theta(t_0 + \varepsilon) < \theta_{Rmax}$ **then**    ▷ Fig. 5(c)
16:         $t_{\theta min}, t_{\theta max} \leftarrow t_0, \infty$
17:       **else if** $\theta_{Rmin} < \theta(t_0 - \varepsilon) < \theta_{Rmax}$ **then**    ▷ Fig. 5(c)
18:         $t_{\theta min}, t_{\theta max} \leftarrow -\infty, t_0$
19:       **else**
20:         $t_{\theta min}, t_{\theta max} \leftarrow$ NaN, NaN
21:       **end if**
22:     **end if**
23:   **else if** one of $t_{Rmin\pm,Rmax\pm}$ valid ▷ Fig. 5(d), (e)
24:     **if** $\theta_{Rmin} < \theta[(\text{valid } t_{Rmin\pm,Rmax\pm}) + \varepsilon] < \theta_{Rmax}$ **then**
      $t_{\theta min}, t_{\theta max} \leftarrow$ valid $t_{Rmin\pm,Rmax\pm}, \infty$
25:     **else**
26:       $t_{\theta min}, t_{\theta max} \leftarrow -\infty$, valid $t_{Rmin\pm,Rmax\pm}$
27:     **end if**
28:   **else if** one of $t_{Rmin\pm}$ and $t_{Rmax\pm}$ valid    ▷ Fig. 5(f)
29:     $t_{\theta min}, t_{\theta max} \leftarrow$ min. and max. of valid $t_{Rmin\pm,Rmax\pm}$
30:   **else if** one of $t_{Rmin\pm,Rmax\pm}$ invalid **then** ▷Fig. 5(g)-(j)
31:     $t_{alone} \leftarrow$ valid solution of $t_{Rmin\pm}$ or $t_{Rmax\pm}$ pair containing invalid solution
32:     $t_{together} \leftarrow$ any solution of $t_{Rmin\pm}$ or $t_{Rmax\pm}$ pair containing two valid solutions
33:     **if** $t_{together} - t_{alone} > 0$ **then**
34:       $t_{\theta min}, t_{\theta max} \leftarrow t_{alone}, \infty$
35:     **else**
36:       $t_{\theta min}, t_{\theta max} \leftarrow -\infty, t_{alone}$
37:     **end if**
38:   **else if** all $t_{Rmin\pm,Rmax\pm}$ valid **then** ▷ Fig. 5(k)-(n)
39:     $t_{\theta min}, t_{\theta max} \leftarrow$ min. and max. of $t_{Rmin\pm,Rmax\pm}$
40:   **else if** all $t_{Rmin\pm}$ and none of $t_{Rmax\pm}$ valid **then**
41:     **if** $\theta_{Rmin} < \pi/2$ **then** ▷ Fig. 5(o), (p)
42:       $t_{\theta min}, t_{\theta max} \leftarrow -\infty, \infty$
43:     **else** ▷ Fig. 5(q)
44:       $t_{\theta min}, t_{\theta max} \leftarrow$ min. and max. of $t_{Rmin\pm}$
45:     **end if**
46:   **else if** all $t_{Rmax\pm}$ and none of $t_{Rmin\pm}$ valid **then**
47:     **if** $\theta_{Rmax} < \pi/2$ **then** ▷ Fig. 5(r)
48:       $t_{\theta min}, t_{\theta max} \leftarrow$ min. and max. of $t_{Rmax\pm}$
49:     **else** ▷ Fig. 5(s), (t)
50:       $t_{\theta min}, t_{\theta max} \leftarrow -\infty, \infty$
51:     **end if**
52:   **else**
53:     $t_{\theta min}, t_{\theta max} \leftarrow$ NaN, NaN
54:   **end if**

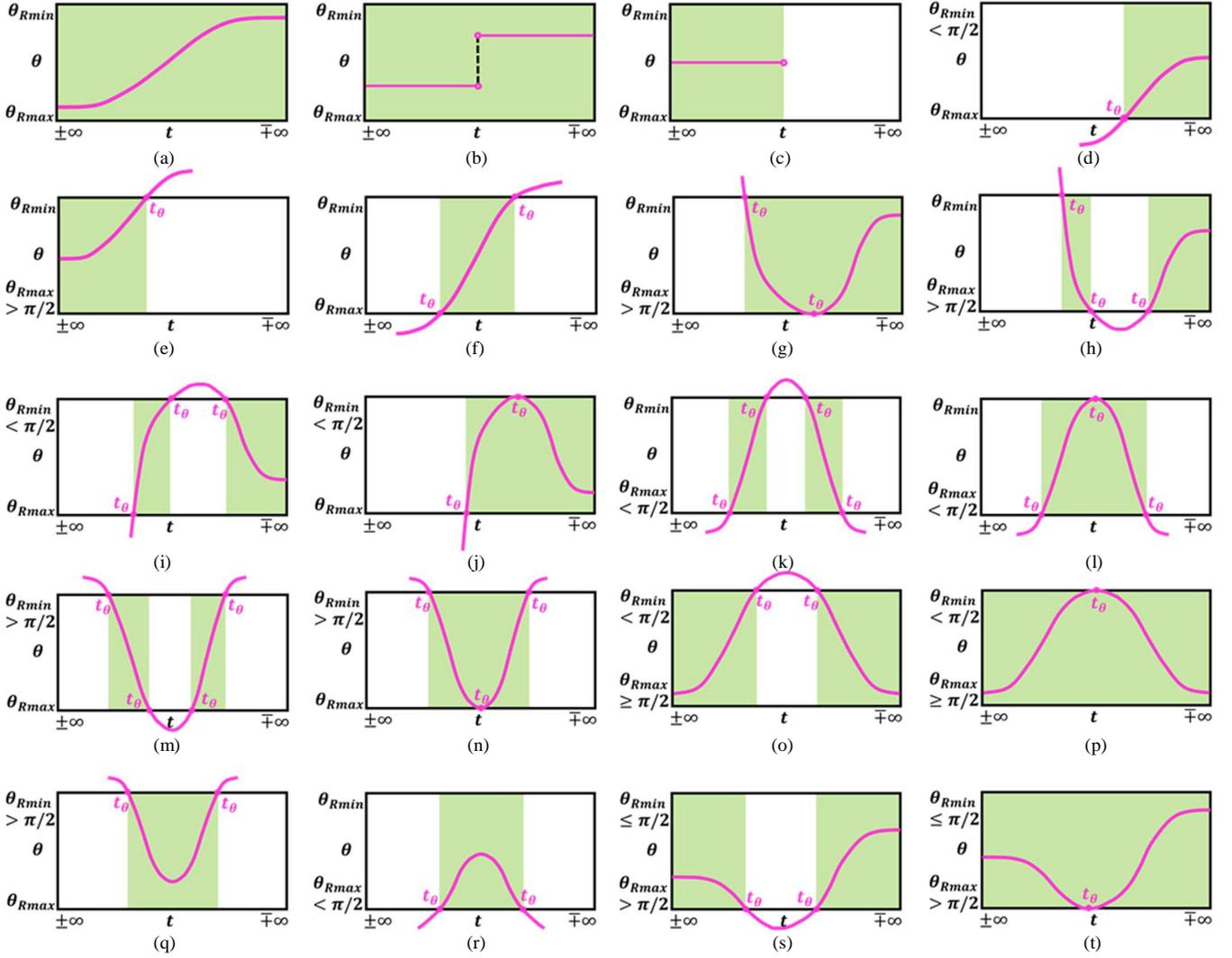

Fig. 5. 20 possible cases where whole or a part of edge line can be located inside the $\theta$ domain of the reflection AZB rectangle. Green part is range of $t$ inside $\theta$ domain of the reflection AZB rectangle.

**Algorithm 2** Continued.
55:   **return** $t_{\theta min}, t_{\theta max}$
56: **end function**

quite challenging and time-consuming process. Instead, we calculate these values for the bounding region containing the facet as shown in Fig. 7. This bounding region is $z'$-directed column-shaped top with two arcs whose center is the diffracting edge. We note that the radiuses of the outer and inner arcs are the horizontally ($x'y'$) furthest and closest distance from the edge to the mesh, respectively, and the minimum and maximum $z'$ components of the bounding region are the minimum and maximum $z'$ components of the facet. We note that $t_{min}$ and $t_{max}$ of the bounding region is determined by four arcs, i.e., lower and upper two arcs. Also, every point on each arc is mapped on the same $t$ value because the diffraction rays are spread in the shape of Keller's cone [27]. To determine the $t$ value where the arc is mapped, we need to solve two equations come from Keller's law [28], that is, $\mathbf{TD} \cdot \mathbf{B} = \mathbf{DR} \cdot \mathbf{B}$ and $\mathbf{D} = \mathbf{P_1} + t\mathbf{B}$ where T, R, D, $P_1$, and $\mathbf{B}$ are locations of Tx, any point on the arc, diffraction point, edge start point, and vector from the edge start point to end point, respectively, and these are shown in Fig.7. After performing some algebra with above equations, we formulate the $t$ where the arc is mapped as follows:

$$t_{arc} = \frac{D_{z'} - P_{1z'}}{P_{2z'} - P_{1z'}} \quad (14)$$

where

$$D_{z'} = \frac{r_T R_{z'} + r_R T_{z'}}{r_T + r_R} \quad (15)$$

and $r_T, r_R$, which are shown in Fig. 7, are horizontal distance between edge and T, R, respectively.

Using (14) and (15), we can calculate $t_{min}$ and $t_{max}$ of the

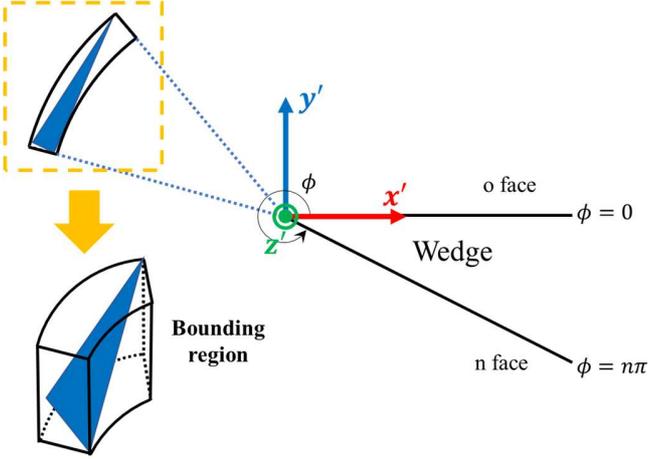

Fig. 6. New basis and bounding region for diffraction-reflection case of ARS.

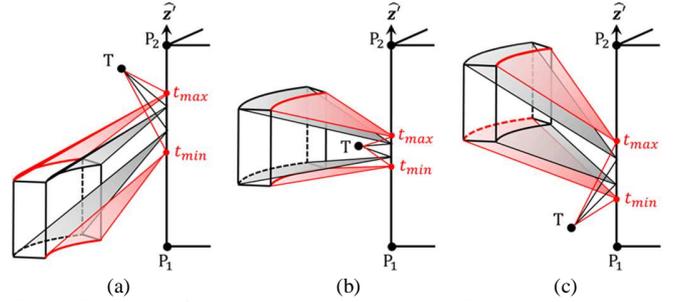

Fig. 8. Three cases for the relative position between Tx and bounding region of diffraction-reflection case of ARS: $T_{z'}$ is (a) larger than maximum, (b) larger than minimum and smaller than maximum, (c) smaller than minimum $z'$ component of the bounding region.

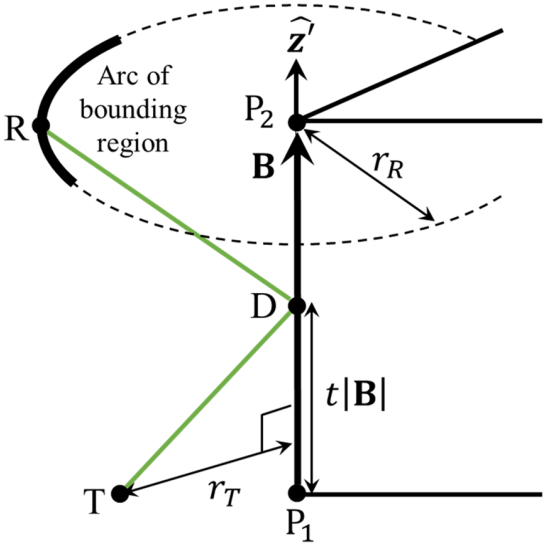

Fig. 7. Diffraction geometry for arc of bounding region of diffraction-reflection case of ARS.

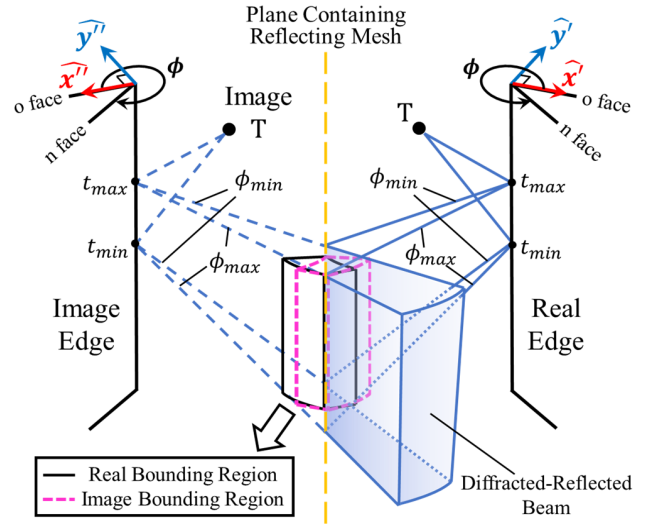

Fig. 9. Image source, edge, bounding region, and new basis for diffraction-reflection case of ARS.

bounding region by calculating all four arcs of the bounding region and comparing their magnitude. However, considering relative position between Tx and bounding region, we can calculate $t_{min}$ and $t_{max}$ directly without magnitude comparison by calculation of $t_{arc}$ for only two arcs. There are three cases for the relative position between Tx and bounding region as shown in Fig. 8. Fig. 8(a), (b), and (c) show where $T_{z'}$ is larger than maximum, larger than minimum and smaller than maximum, and smaller than minimum $z'$ component of the bounding region. For each case, $t_{min}$ and $t_{max}$ value occurs at the lower inner and upper outside, lower and upper inside, and lower outside and upper inside arcs, respectively. Eventually, the proposed bounding region shape has an advantage of simple calculation of $t_{min}$ and $t_{max}$ value using (14) and (15) considering only two arcs of the bounding region. This is a superior characteristic compared to using the popular axis-aligned bounding box which has eight corners and whose corner where the $t_{min}$ and $t_{max}$ values are mapped is difficult to predict and the formulation of the $t_{min}$ and $t_{max}$ values of each corner is more complex than (14) and (15).

After determining the first and second $\phi - t$ AZB rectangles, we determine the intersection between the two. Subsequently, we generate images of the source and edge in the reflecting facet and transform basis similarly with the multiple reflection case, i.e., the first and second new basis $\mathbf{x}''$ and $\mathbf{y}''$ are $\mathbf{x}'$ and $\mathbf{y}'$ mirrored to the reflecting facet, respectively, as shown in Fig. 9. Then, the problem can be reduced to a single diffraction of the image source and image edge on the image bounding region, i.e., we can directly use the intersected AZB rectangle, that is, shrunk AZB rectangle with margins of $\phi_{min}, \phi_{max}, t_{min}$, and $t_{max}$, as shown in Fig. 9.

For the diffraction-multiple reflection, a same procedure can be used to shrink the AZB rectangle because every diffraction-multiple reflection problem can be reduced to the single diffraction of image source and image edge as illustrated above. The only difference is that the two AZB rectangles (to calculate the intersection) for each bouncing order are the shrunk AZB rectangle of former bouncing order and one for the reflecting facet and the image edge. In this way, size of the AZB rectangle can be shrunk exponentially as the bouncing order increases.

In this section, the ARS method for the multiple reflection,

reflection-diffraction, and diffraction-reflection is described. The ARS method of other propagation sequences, e.g., multiple reflection-diffraction, reflection-diffraction-multiple reflection, and etc., can be formed by combinations of the above three propagation sequences.

### III. HETEROGENEOUS COMPUTING BASED ARS-ACCELERATED IMAGE THEORY (HAIT) METHOD RAY TRACING FRAMEWORK

In this section, we describe our novel HAIT method ray tracing framework. The framework is divided into three part; 1. Visibility preprocessing. 2. Visibility tree generation using ARS and CPU parallel computing. 3. ST and field calculation using CPU/GPU heterogeneous computing.

*A. Visibility preprocessing*

The first part of the HAIT is the visibility preprocessing. In this part, the framework derives the visibility relationships between facet-facet, facet-edge, Tx-facet, Tx-edge, facet-FOP, and edge-FOP. The preceding four and following two relationships are used to accelerate the visibility tree generation and the ST, respectively. In this paper, we do not describe our visibility preprocessing algorithm in detail considering the size of the paper. However, readers can refer some visibility preprocessing algorithms in [25] and [29]-[33].

*B. Visibility tree generation using ARS*

The second part of the HAIT is the visibility tree generation. In this part, we generate the visibility tree containing every possible sequence of primitives of the geometry object, i.e., facet or edge, that ray can bounce as shown in Fig. 10. The level 1 of the tree is directly constructed using Tx-facet and Tx-edge visibility relationship obtained from the visibility preprocessing. The higher levels of the tree can be efficiently built using the ARS. Furthermore, facet-facet and facet-edge visibility relationships obtained from the visibility preprocessing accelerate the visibility tree generation one more time by considering only visible facets and edges from the facet or edge that produce beam when we find facets or edges inside the beam.

Also, we use Open Multi-Processing (OpenMP), which is a CPU parallel computing application programming interface, for parallel computing acceleration. Because every subtree of the root node of the tree is independent with each other, we can parallelize the tree generation operation by let each OpenMP thread build different subtrees simultaneously, resulting in maximally n times faster computation than general serial computation when n CPU threads are used.

We note that we do not generate DAZB matrices which are used by general AZB method to reduce the number of ray-facet intersection test. It is because DAZB matrices should be generated for every image Tx and it results into usage of prohibitively large amount of memory; therefore, it is not suitable for the parallel computing.

We also note that when we build the visibility tree, we do not attach FOPs to the tree as shown in Fig. 10. It is because if we attach FOPs to the tree during the tree construction, the tree size would be extremely large when an analysis scenario treats large amount of field observation point, resulting in requirement of prohibitively large amount of memory and precluding the parallel computing acceleration. Also, finding thousands of FOPs inside the AZB rectangle is a simple but assembly of greatly large amount identical operations and it is more suitable for GPU parallel computing with thousands of cores rather than CPU parallel computing. Therefore, we attach FOPs bit by bit to the visibility tree in the following ST part and conduct ST using GPU parallel computing. This trick that ignoring FOPs in the tree generation part has an additional advantage that the completed visibility tree without FOPs can be reused when the locations of FOPs are changed because the tree is independent with the location of FOPs. It is a superior function that the SBR method RT does not retain because the general SBR method does not construct the visibility tree.

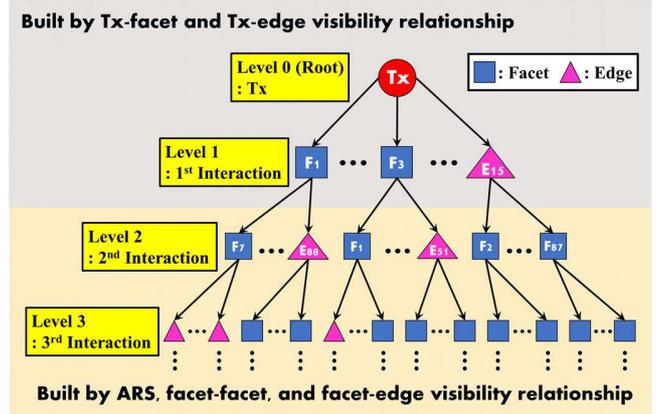

Fig. 10. Visibility tree without field observation points.

*C. ST and field calculation*

The final part of HAIT is ST and field calculation part. In this part, we propose a GPU and CPU heterogeneous computing where the GPU and CPU simultaneously conduct ST and field calculation, respectively.

First of all, because HAIT does not generate DAZB matrices, we need to use other acceleration structure, e.g., kD-tree and bounding volume hierarchy (BVH) structure. We adopted the BVH structure which is provided by NVIDIA® OptiX™. OptiX is a RT engine that is a programmable system designed for NVIDIA GPUs and other highly parallel architectures. It provides high performance BVH structure, ray-facet intersection test algorithm, and high-quality optimization at low-level. Also, new versions of OptiX are periodically released, so that periodical performance improvement is possible by small change of program code. Other detailed information of OptiX can be found in [34]. Using OptiX, we conduct ST in the GPU environment. Detailed ST algorithm is illustrated in [13].

For field calculation, we use the geometrical optics (GO) [24], [28], [35]-[39] and uniform geometrical theory of diffraction (UTD) [40], [41] to calculate reflection and diffraction fields, respectively.

Fig. 11 shows the flowchart of the ST and field calculation part. First of all, we create BVH structure using OptiX to accelerate the ST. Next, we attach FOPs to every node of the visibility tree. For the root node, we attach every FOPs, and for the other nodes in higher levels, we attach FOPs visible from facet or edge corresponding each node, considering facet-FOP, edge-FOP visibility relationship obtained from the visibility preprocessing. The maximum number of FOP to be attached at once is $2^{30}$ because the maximum ray launch size of OptiX is $2^{30}$ [42]. It means that if the maximum number of required ST is over $2^{30}$, we cannot perform ST at once and should partition the algorithm into multiple ST sets. On the basis of this fact, we propose a heterogeneous computing algorithm. This heterogeneous computing simultaneously performs $n^{th}$ ST set using OptiX, i.e., GPU parallel computing, and field calculation for ST pass ray paths of $(n-1)^{th}$ ST set using OpenMP, i.e., CPU parallel computing, and each GPU/CPU thread performs ST/field calculation for each branch of the tree. This is possible because the compute unified device architecture (CUDA) returns the control to the CPU right after invoking CUDA kernel without waiting for the completion of GPU computation [43], resulting in increase of computation efficiency.

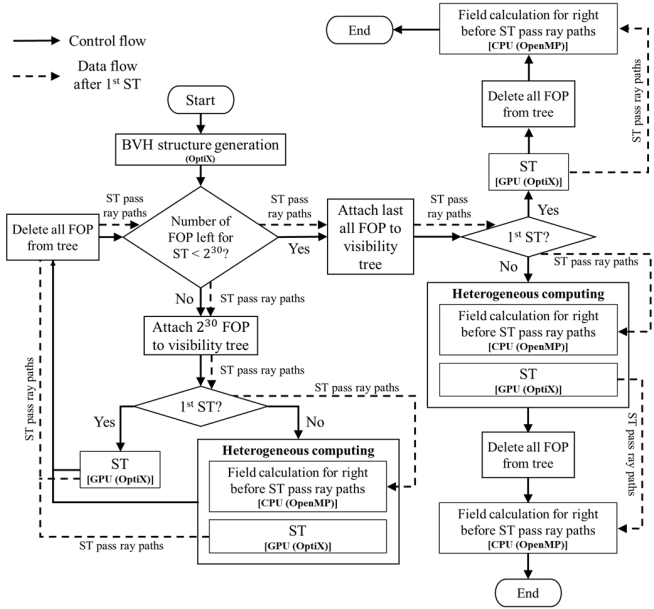

Fig. 11. Flowchart of the ST and field calculation part.

We note that FOP attachment/ST/field computation is started from the node of low level to the high level of the visibility tree to escape CUDA thread divergence as much as possible [44].

## IV. VALIDATION

In this section, we describe the performance of the proposed HAIT by comparison between the HAIT, HAIT without ARS, and an IT solver of commercial ray tracer WinProp about two outdoor scenarios: simple and complex ones. All simulations were conducted on an identical computation environment: Intel(R) Xeon(R)CPU E5-2687 v4 @ 3.00 GHz 3.00 GHz (2 processor) 512GB RAM (16 cores are used) and NVIDIA RTX A6000 GPU. The HAIT and HAIT without ARS are implemented in C/C++ language and OptiX.

The first simulation scenario is a simple scenario which is shown in Fig. 12 and Table I. Fig. 12 shows the simulation environment describing Teheran-ro in Seoul, South Korea. It is constructed by AutoCAD and contains 1,218 triangle facets. In this scenario, E-field is analyzed at 28 GHz for a Hertzian dipole Tx radiating power of 1 W. Tx is located at (-16.6, -94.3, 10) m and 100 FOPs are uniformly distributed on the line between (-143, -81, 1.5) and (242, 51, 1.5) m. Also, the material of the simulation environment is set to perfect electric conductor (PEC) and it is because, for the diffraction of wedges with finite conductivity, our simulator and WinProp uses different heuristic diffraction coefficient models, i.e., models in [41] and [45] respectively, resulting in difference simulation results. Lastly, 1 to 6 of maximum number of reflections, combined with one order of diffraction are considered. Table II shows the performance of HAIT, HAIT without ARS, and WinProp for the simple scenario. Naturally, the HAIT and HAIT without ARS shows zero root mean square percentage error (RMSPE) because the ARS removes only unnecessary facets and edges from the visibility tree. Although HAIT shows similar computation time with the HAIT without ARS at the maximum bouncing order up to 4, the HAIT shows about 1.6 and 3 times faster computation speed than the HAIT without ARS at the maximum bouncing order of 5 and 6, respectively. This phenomenon is easily expected because, as described in Section. II, the effect of the ARS method increases exponentially as the bouncing order increases. Looking at the simulation results of WinProp, we note that the simulation is conducted only for the maximum bouncing orders up to 3 because it consumed tremendous computation time at the maximum bouncing order over 3. The RMSPEs between HAIT and WinProp are very small, i.e., under 1.5 %, and this small error is probably because of the usage of different data types in each ray tracer, resulting in floating-point number precision error. Fig. 13 shows the E-field simulation results of the HAIT, HAIT without ARS, and WinProp for the maximum bouncing order of 3. Although the HAIT shows about two times faster computation time than WinProp at the maximum bouncing order of 1 and 2, about 651 times faster computation time is confirmed at the maximum bouncing number of 3. The reason of this huge difference is hard to explain because the exact algorithm of WinProp is not open to the public. However, authors expect that it is because the IT solver of WinProp is a general IT ray tracer, and like any other general IT ray tracer, its computation time exponentially increases at the high maximum bouncing order and large amount of the FOPs. On the other hand, the HAIT uses the ARS method and heterogeneous computing-based ST and field calculation algorithm optimized for large number of FOPs, resulting in the superior performance by large amount reduction and acceleration of ST. Resultingly, authors expect that the ARS and heterogeneous computing algorithm creates the far superior performance compared to WinProp.

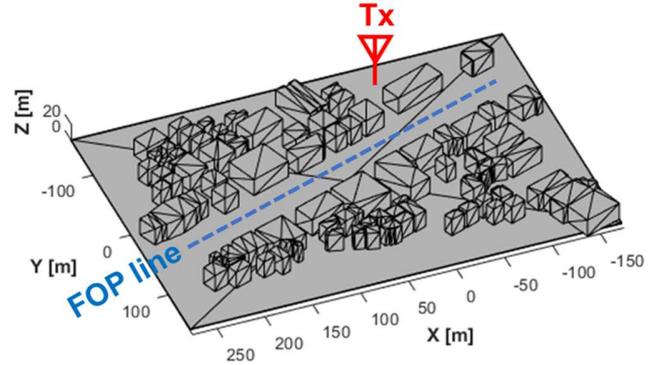

Fig. 12. Simulation environment of the first, i.e., simple, scenario describing Teheran-ro in Seoul, South Korea and locations of Tx and FOPs of simple scenario.

TABLE I
SIMULATION PARAMETERS OF SIMPLE SCENARIO

| Frequency | Tx antenna | Tx radiation power |
|---|---|---|
| 28 [GHz] | Hertzian dipole (z-directed) | 1 [W] |
| Max. Bouncing # | Tx location | Number of FOP |
| 1, 2, 3, 4, 5, 6 | (-16.6, -94.3, 10) [m] | 100 |
| Location of FOP | | |
| Uniformly distributed on line between (-143, -81, 1.5) and (242, 51, 1.5) [m] | | |

TABLE II
PERFORMANCE FOR SIMPLE SCENARIO

| Max. bouncing # | RMSPE [%] | | Computation time [sec] | | |
|---|---|---|---|---|---|
| | HAIT vs HAIT w/o ARS | HAIT vs WinProp | HAIT | HAIT w/o ARS | WinProp |
| 1 | 0 | 0.99 | 10 | 10 | 19 |
| 2 | 0 | 0.81 | 10 | 10 | 20 |
| 3 | 0 | 1.09 | 10 | 11 | 6,512 |
| 4 | 0 | - | 12 | 13 | - |
| 5 | 0 | - | 17 | 28 | - |
| 6 | 0 | - | 52 | 158 | - |

Fig. 14 shows the simulation environment of the second, i.e., complex, scenario describing Gangnam in Seoul, South Korea. It contains 4,996 triangle facets and is broader, i.e., 1 km x 1 km, and more complex urban area than the simple scenario. Table III shows the simulation parameters. In this scenario, similar to the simple scenario, E-field is analyzed at 28 GHz for a Hertzian dipole Tx radiating power of 1 W. Tx is located at (-46, 14, 10) m and 10,000 FOPs are uniformly distributed on the area (x: -560~440, y: -430~570, z: 1.5) m. Also, due to the same reason in the simple scenario, the material of the simulation environment is set to PEC. Lastly, 1 to 6 of maximum order of reflections, combined with one order of diffraction are considered. Table IV shows the performance of the complex scenario. Similar to the simple scenario, the HAIT and HAIT without ARS shows zero RMSPE. Although HAIT shows similar computation time and memory consumption with the HAIT without ARS at the maximum bouncing order up to 3, the HAIT shows about 1.5 times faster computation speed and 1.4, 1.2 GB less memory consumption of CPU, GPU than the HAIT without ARS at the maximum bouncing order of 4. Moreover, at the maximum bouncing order of 5, the memory consumption of HAIT without ARS is exponentially increased and encounters out of memory error at the GPU environment, i.e.

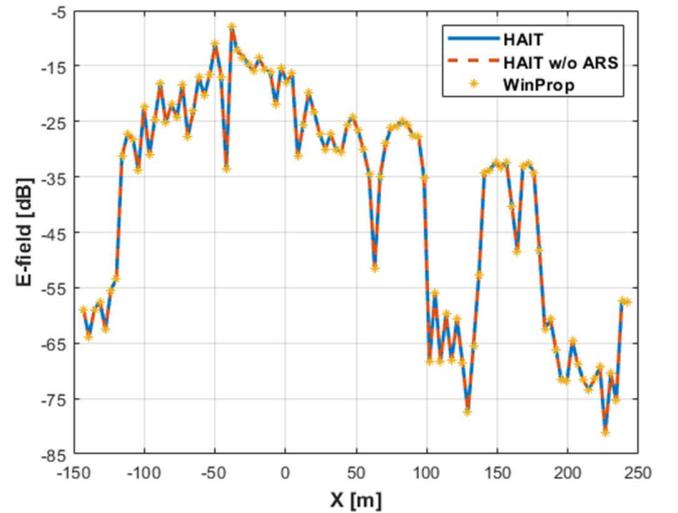

Fig. 13. E-field simulation results of simple scenario with maximum ray bouncing order of 3 for HAIT, HAIT without ARS, and WinProp .

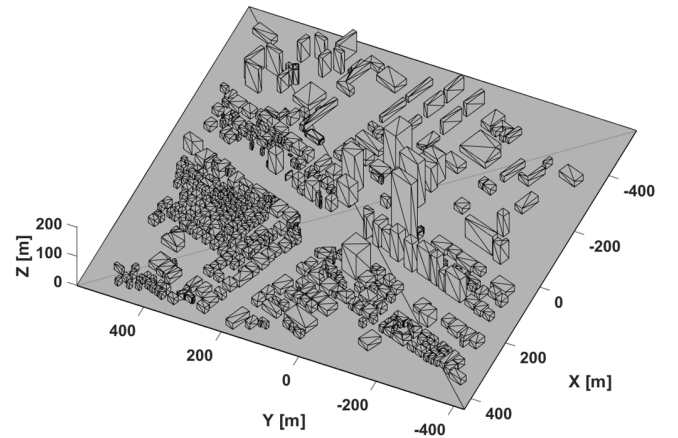

Fig. 14. Simulation environment of the second, i.e., complex, scenario describing Gangnam in Seoul, South Korea.

TABLE III
SIMULATION PARAMETERS OF COMPLEX SCENARIO

| Frequency | Tx antenna | Tx radiation power |
|---|---|---|
| 28 [GHz] | Hertzian dipole (z-directed) | 1 [W] |
| Max. Bouncing # | Tx location | Number of FOP |
| 1, 2, 3, 4, 5, 6 | (46,14,10) [m] | 10,000 |
| Location of FOP | | |
| Uniformly distributed on area (x: -560~440, y: -430~570, z: 1.5) [m] | | |

GPU memory consumption is larger than 48GB, while the HAIT consumes only 19.8GB at the GPU environment. At the maximum ray bouncing order 6, the HAIT also encounters GPU out of memory error. To manage this error, authors conducted some serialization by partitioning the visibility tree into 6 parts. More specifically, in the complex scenario, the root node of the full visibility tree has 513 child nodes. Authors firstly completed 85 subtrees of the root node and conducted ST and field calculation for these 85 subtrees. After that, authors conducted other 85 or 86 subtrees construction and ST and field calculation sequentially for 5 more times to deal with every subtrees of the root node. This visibility tree partitioning trick can be used to control the trade-off between the memory consumption and simulation time. Fig. 15 shows the E-field simulation results of the HAIT with maximum bouncing number of 6. To the best of authors knowledge, it is the first IT ray tracer that is feasible for the simulation of the such a large-scale, i.e., 1km x 1km wide, and dense outdoor environment with high maximum ray bouncing order of 6 and the large amount of field observation points of 10,000. Looking at the simulation results of WinProp, we note that the simulation is conducted only for the maximum bouncing orders up to 2 because it consumed tremendous computation time at the maximum bouncing order over 2. The RMSPEs between HAIT and WinProp are very small, i.e., under 4 %, and this small error is probably because of the floating-point number precision error as described in the simple scenario. Although the HAIT shows slower computation time than WinProp at the maximum bouncing order of 1 because of the visibility preprocessing process of HAIT which takes 145 seconds, it is confirmed that HAIT is faster than WinProp about 34 times at the maximum bouncing number of 2.

TABLE IV
PERFORMANCE OF COMPLEX SCENARIO

| Max. Bounce # | RMSPE [%] | | Computation time [sec] | | | Memory [GB] | | | | | |
| --- | --- | --- | --- | --- | --- | --- | --- | --- | --- | --- | --- |
| | HAIT vs HAIT w/o ARS | HAIT Vs WinProp | HAIT | HAIT w/o ARS | WinProp | HAIT | | HAIT w/o ARS | | WinProp | |
| | | | | | | CPU | GPU | CPU | GPU | CPU | GPU |
| 1 | 0 | 3.95 | 154 | | 28 | 0.73 | 1.22 | 0.73 | 1.22 | 0.14 | - |
| 2 | 0 | 3.33 | 155 | 155 | 5,281 | 0.75 | 1.51 | 0.75 | 1.51 | 0.24 | - |
| 3 | 0 | - | 188 | 194 | - | 21.4 | 14.3 | 21.4 | 14.3 | - | - |
| 4 | 0 | - | 508 | 788 | - | 25.7 | 14.7 | 27.1 | 15.9 | - | - |
| 5 | - | - | 5,461 | - | - | 31.3 | 19.8 | - | > 48 | - | - |
| 6 | - | - | (Tree partition) 73,688 | - | - | (Tree partition) 49.8 | (Tree partition) 39.9 | - | - | - | - |

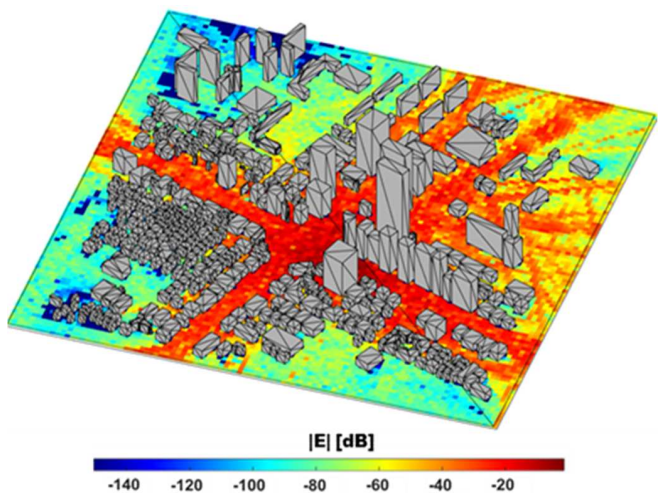

Fig. 15. E-field simulation results of the HAIT with maximum bouncing number of 6 for complex scenario.

## V. CONCLUSION

In this paper, the HAIT ray tracing framework for the outdoor propagation modeling is proposed. The HAIT is featured by ARS method and the CPU/GPU heterogeneous computing. The ARS method effectively accelerates and reduces the generation time and size of the visibility tree by letting the beam radiate only to the portion that is illuminated by the beam of the former bouncing order. Resultingly, compared to the case without ARS method, more than 3 times faster simulation time is confirmed in the simple scenario. Also, it is confirmed that, using ARS, it is possible to handle the complex scenario, that is impossible for without ARS case due to the out of memory error, by using lesser than half of memory consumption. This efficiency improvement increases exponentially further as the maximum ray bouncing order, complexity of simulation geometry, and the number of FOPs increase.

The CPU/GPU heterogeneous computing accelerates overall part of algorithm by letting OpenMP CPU parallel computing handle the visibility tree generation and field computation task, OptiX GPU parallel computing handle the ST. At this time, it was able to utilize GPU parallel computing and handle thousands of FOPs using following two tricks: 1. Instead of creating DAZB matrices while the visibility tree generation process, we use BVH acceleration structure in the ST part. 2. We ignore the FOPs at the visibility tree generation part and attach them at the ST part. Eventually, it is confirmed that the HAIT is more than 651 times faster than the IT solver of WinProp.

Until now, despite their high accuracy, the utilization of the conventional IT ray tracers was limited to simple simulation environments with small number of FOPs and low maximum ray bouncing order due to their poor computational efficiency. However, in this paper, it is confirmed that the proposed HAIT can handle 1 km x 1km wide complex urban environment with the high maximum ray bouncing order of 6 and thousands of FOPs. Therefore, authors expect that it could be used without doubt about systemic error to areas that was exclusive to the SBR method ray tracers due to their demand of complex analysis, such as channel modeling and coverage analysis. The frequency and computing power we use are continuously increasing these days and IT ray tracers, therefore, would be one of the good choices for the complex outdoor propagation modeling. The HAIT would be a cornerstone for it. There are two promising future works. First one is applying additional propagation mechanisms such as transmission and multiple diffraction. Second one is applying additional computing techniques, such as cluster computing to ARS part, and multiple GPU parallel computing to ST/field calculation part, for further computation acceleration.

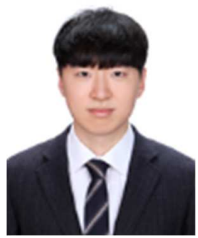

**Yongwan Kim** received a B.S. degree in electrical engineering from Chungnam National University, Korea, in 2019. He is currently pursuing a Ph.D. degree with Seoul National University, Korea. His current research interests include ray-tracing EM analysis technique, EM theory, parallel computing, heterogeneous computing, and EMI/EMC.

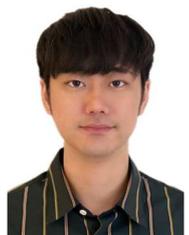

**Hyunjun Yang** received a B.S. degree in electronic and electrical engineering from Pohang University of Science and Technology (POSTECH), Pohang, South Korea, in 2020. He is currently pursuing an integrated master's and Ph.D. degree in the Department of Electrical Engineering and Computer Science at Seoul National University, Seoul, South Korea. His current research interests include Ray tracing EM analysis technique, EM Theory, wireless propagation channel measurement & modeling and RIS.

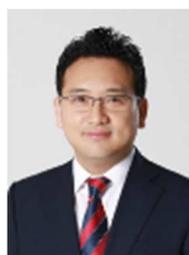

**Jungsuek Oh** (S'08) received his B.S. and M.S. degrees from Seoul National University, Korea, in 2002 and 2007, respectively, and a Ph.D. degree from the University of Michigan at Ann Arbor in 2012. From 2007 to 2008, he was with Korea Telecom as a hardware research engineer, working on the development of flexible RF devices. In 2012, he was a postdoctoral research fellow in the Radiation Laboratory at the University of Michigan. From 2013 to 2014, he was a staff RF engineer with Samsung Research America, Dallas, working as a project leader for the 5G/millimeter-wave antenna system. From 2015 to 2018, he was a faculty member in the Department of Electronic Engineering at Inha University in South Korea. He is currently an Associate Professor in the School of Electrical and Computer Engineering, Seoul National University, South Korea. His research areas include mmWave beam focusing/shaping techniques, antenna miniaturization for integrated systems, and radio propagation modeling for indoor scenarios. He is the recipient of the 2011 Rackham Predoctoral Fellowship Award at the University of Michigan. He has published more than 50 technical journal and conference papers, and has served as a technical reviewer for the IEEE Transactions on Antennas and Propagation, IEEE Antenna and Wireless Propagation Letters, and so on. He has served as a TPC member and as a session chair for the IEEE AP-S/USNC-URSI and ISAP. He has been a senior member of IEEE since 2017.